\documentclass{siamltex}

\usepackage[T1]{fontenc}
\usepackage{latexsym}
\usepackage{graphicx}
\usepackage{beton}
\usepackage{euler}
\usepackage{stmaryrd}
\usepackage{harvard}
\usepackage{amsmath,amssymb}

\newcommand{\R}{{\mathbb R}}
\newcommand{\G}{{\mathbb F}}

\newcommand{\meas}[1]{\left\llbracket#1\right\rrbracket}
\newcommand{\norm}[1]{\left\|#1\right\|}
\newcommand{\xmapsto}[1]{\stackrel{#1}{\longmapsto}}

\makeatletter
\newcommand{\cond}{{\operator@font cond}}
\makeatother

\renewcommand{\i}{{\rm i}}

\title{A Model for Understanding Numerical Stability}

\author{Folkmar Bornemann\thanks{Zentrum Mathematik, Technische Universit\"{a}t M\"{u}nchen,
        Boltzmannstr. 3, 85747 Garching, Germany ({\tt bornemann@ma.tum.de}). \emph{e-print arXiv:math.NA/0503100
        as of \today}.}}

\begin{document}

\maketitle

\begin{abstract}
We present a  model of roundoff error analysis that combines simplicity with predictive power. Though
not considering all sources of roundoff within an algorithm, the model is related to a recursive
roundoff error analysis and therefore capable of correctly predicting stability or instability of an algorithm.
By means of nontrivial examples, such as the componentwise backward stability
analysis of Gaussian elimination with a single iterative refinement step, we demonstrate that the model even yields
 quantitative backward error bounds that show all the known problem-dependent terms (with the exception of
dimension-dependent constants, which are the weak spot of any a priori analysis). The model can serve as a
convenient tool for teaching or as a heuristic device to discover stability results before entering a further, detailed analysis.
\end{abstract}

\section{Introduction}

An algorithm for the numerical evaluation of a complicated function $f$
is just a  decomposition into simple intermediate steps, such as arithmetic
operations, elementary transcendental functions, or well-behaved and well-understood library algorithms (e.g., matrix multiplication):
\[
f = g_1 \circ g_2 \circ \cdots \circ g_k.
\]
In the realm of floating point arithmetic, each of these intermediate steps is contaminated by
roundoff and hence contributes to the final perturbation of the result in a twofold fashion:
first, by \emph{generating} roundoff error itself, and second, by \emph{propagating} the roundoff errors of previous steps.
Since the early days of numerical computing there has been much progress in clarifying the underlying structure and organizing the results in a concise,
easily interpreted form.
However, a detailed analysis \cite{Higham} is still often quite involved and remains a battle-field for experts, too tedious to teach
and explain in detail beyond the most trivial cases in a beginner's course on numerical analysis. The instructor typically chooses
between two options: skipping the nontrivial results (such as stability of Gaussian elimination) entirely, or just stating the results
without proof. Either choice is unsatisfactory for good students since they cannot develop an {\em understanding} of the mathematical structure and reasons.

We will demonstrate in this paper, that the overall behavior of an algorithm can very often be well understood
by analyzing a simplified {\em model} of the sources of roundoff error. As in the natural sciences
such a model has to balance simplicity with predictive power. If such a simple model basically
leads to the same predictions, qualitatively and perhaps even quantitatively, as a full-fledged
a priori roundoff error analysis, we may rightly claim to have contributed to the understanding
of the algorithm's behavior. In fact, all the estimates of our model analysis that we present in this paper
will give the \emph{same} estimates as a detailed a priori analysis---with the only exception of the dimension-dependent
constants, which are, however, anyway the weak spot, and therefore least important part, of any roundoff error analysis \cite[p.~65]{Higham}.
In particular, with just a few lines of simple
calculations we will obtain the nontrivial results on the norm- and componentwise backward stability of Gaussian elimination
ranging from the early work \cite{Wilkinson} to the analysis of iterative refinement \cite{Skeel2}.

In addition to being a convenient (and to the experience of the author also successful) tool in teaching, our model might
serve as a \emph{heuristic} device in {discovering} the structure of a stability result---before one enters, in a second step, taking advantage of the obtained knowledge,
a fully detailed roundoff error analysis.

\paragraph{The Model}
The roundoff error analysis that we propose is based on the observation that in many if not most cases
a critical intermediate step can be identified that leads to a natural decomposition
\[
f = \underbrace{g_1 \circ \cdots \circ g_j}_{=h} \circ \underbrace{g_{j+1} \circ \cdots \circ g_k}
_{=g} = h \circ g
\]
into just two fundamental steps. Now, the model is based on the \emph{simplifying assumption} that roundoff error just affects
the single intermediate result---after being output by $g$, before being input to $h$. That is, we analyze
the error of the \emph{realization map}
\[
\tilde{f} = h \circ fl \circ g.
\]
Here $fl: \R^p \to \G^p$ denotes componentwise rounding, subject to the standard model
of floating point arithmetic
\[
| fl(x) -x | \leq u\cdot |x|,
\]
where $u$ denotes the unit roundoff of the arithmetic ($u\approx 1.11\times 10^{-16}$ for
IEEE double precision) and $\G$ the floating point numbers. We understand $|x|$ componentwise
for vectors and matrices.

\paragraph{Outline of the Paper}
In {\S}\ref{sect.backward} we analyze the {backward} stability of the realization map~$\tilde f$, which
turns out to be determined by the condition number of $g^{-1}$. We will specify the relation of the model to a complete analysis. In fact,
if the model is unstable the same has to be expected for the real situation. On the other hand, if $g$ and $h$ are realized
by backward stable algorithms, the resulting algorithm for $f$ would inherit the stability of the model. This helps to understand the
success of our model and suggests a recursive approach to a full roundoff error analysis.

The rest of the paper studies some algorithms for the solution of a linear system $Ax=b$. In {\S}\ref{sect:oettli} we recall some classic expressions
for the backward error of linear systems that are the point of departure for the simple estimates to follow.
In {\S}\ref{sect:naiv} we study the na\"{\i}ve algorithm,
that is, multiplication with $A^{-1}$, and show its instability for badly conditioned matrices. In {\S}\ref{sect:GEnorm} we study the normwise
backward error of Gaussian elimination and obtain the classic result \cite{Wilkinson}. In {\S}\ref{sect:GEcomp} we get the result
\cite{Skeel1} on the componentwise backward error of Gaussian elimination, correctly predicting the influence of the scaling of the system.
Finally, in {\S}\ref{sect:GESkeel} we show how to discover within the frame of our model the result \cite{Skeel2} that a single step of iterative refinement implies componentwise
backward stability of Gaussian elimination.

\section{Backward Stability}\label{sect.backward} The main result of a \emph{qualitative} study of our model can be summarized
as follows:
\medskip
\begin{quote}
\emph{Backwards stability requires that $g^{-1}$ is well conditioned.}
\end{quote}
\medskip
In fact, backward stability analysis requires the result of the algorithm for
input
$x$, that is $\tilde{f}(x)$ in our model, to be represented as the \emph{exact} solution to perturbed data:
$\tilde{f}(x) = f(x+\Delta x)$.
Writing $w=g(x)$ for short, we have
\[
\tilde{f}(x) = h(fl(w)) = h(w+ \Delta w),\qquad |\Delta w| = |fl(w) - w| \leq u \cdot |w|.
\]
Assuming $g$ to be invertible, we propagate $\Delta w$ backwards to obtain an estimate for $\Delta x$:
\[
x + \Delta x = g^{-1}(w + \Delta w), \qquad |\Delta x | \leq \kappa_{g^{-1}}\,u\cdot |x| + O(u^2),
\]
where $\kappa_{g^{-1}}$ denotes the (componentwise) relative condition number of $g^{-1}$ at $w$. Hence,
the backbard error is bounded by the unit roundoff amplified by $\kappa_{g^{-1}}$.

\subsection{Examples}\label{sect2:examples}

\paragraph{A} Consider the evaluation of $f(x)=\log^2(1+x)$ for $x\approx 0$. A direct implementation
of the defining formula corresponds to the decomposition
\[
f : x  \xmapsto{g} w = 1+x \xmapsto{h} \log^2(w).
\]
Now, because $w=1+x \approx 1$ the inverse function $g^{-1} : w \mapsto x = w-1$ is a subtraction
in the cancelation regime, thus badly conditioned. Hence, we predict instability of the formula,
which simple examples confirm. In fact,  the bad conditioning of $g^{-1}$ reflects the
{\em loss of information} in $g$: we have $fl(g(x)) = fl(1+x)=1$ as soon as $x$ is smaller
than the resolution of the machine arithmetic.  In general, well-conditioning of $g^{-1}$, however, requires that
the input $x$ is accurately reconstructable from the intermediate result $w=g(x)$.

\paragraph{B} The solution $x \in \R^m$ of a linear system of equations $Ax =b$
with a nonsingular $A \in \R^{m\times m}$ can formally be respresented as $x =A^ {-1}\cdot b$.
This suggest the na\"{\i}ve algorithm corresponding to the decomposition
\[
f : A  \xmapsto{\;\, g \;\,} A^{-1}  \xmapsto{\;\, h \;\,} x = A^{-1}\cdot b.
\]
Now, $g^{-1}: A^{-1} \mapsto A$ is just $g$ again, its condition is (in the normwise case)
the condition number of the matrix $A$. Thus, we expect the algorithm to be
unstable for certain badly conditioned matrices. Examples that display such instability will be given
 in {\S}\ref{sect:naiv} where we extend our analysis to a more quantitative setting.

\paragraph{C} On the other hand, the solution of the linear system $Ax =b$ by Gaussian
elimination  corresponds to the decomposition
\[
f: A\xmapsto{\;\, g \;\,} (L,U)  \xmapsto{\;\, h \;\,} x .
\]
Here, $g$ represents the $LU$-factorization step, whereas
$h$ represents the substitution steps. Now, the inverse of $g$, that is
\[
g^{-1}: (L,U) \longmapsto A = L\cdot U,
\]
is just {\em matrix multiplication}. Its condition number can be estimated by
\[
\kappa_{g^{-1}} \leq 2\, \frac{\|\,|L|\cdot|U|\,\|}{\|A\|},
\]
which is,  as will be discussed in more detail in {\S}\ref{sect:GEnorm}, sufficient to explain the instabilities to be observed for Gaussian elimination
with or without partial pivoting.

\subsection{Relation of the Model to a Complete Analysis}
In fact, the condition number of $g^{-1}$ turns out to be relevant for a full roundoff error analysis, too.
Here, we would recursively define the realization of $f = h \circ g$ by
\[
\tilde{f} = \tilde{h} \circ \tilde{g},
\]
starting with the backward stable realization of the arithmetic operations and basic elementary
functions. (Of course, in general we cannot assume that in each step of this recursion the $g$-part of the
decomposition is invertible. However, it is possible to give a reasonable definition of $\kappa_{g^{-1}}$ even if $g$ is many-to-one.)

With $\meas{\Delta x}$ denoting maximum componentwise relative error\footnote{That is, for the perturbation
 $\Delta x \in \R^m$ of a quantity
$x \in \R^m$ we have $\meas{\Delta x} = \max_{j=1:m} |\Delta x_j|/|x_j|$ with the convention that $0/0=0$.} we define the smallest number $\beta_f\geq 0$ such that
\[
\tilde{f}(x) = f(x+\Delta x),\qquad \meas{\Delta x} \leq \beta_f \cdot u  + O(u^2)
\]
as the \emph{stability indicator} of $\tilde{f}$. Backwards stability requires  $\beta_f$ to be not too large.

\begin{lemma}
For $g$ invertible there holds the recursive estimate
\begin{equation}\label{eq:calculus}
\beta_f \leq \beta_g + \kappa_{g^{-1}}\cdot \beta_h.
\end{equation}\vspace*{-0.5cm}
\end{lemma}
\begin{proof}
The stability indicator of $\tilde{h}$ gives
\[
\tilde{f}(x) = \tilde{h}(\tilde{g}(x)) = h(\tilde{g}(x) + \Delta w),\qquad \meas{\Delta w}\leq \beta_h \cdot u + O(u^2).
\]
The stability indicator of $\tilde{g}$ and the relative condition number of $g^{-1}$ allow for the estimates
\begin{align*}
\tilde{g}(x) + \Delta w &= g(x+\Delta x_1) + \Delta w, \qquad \meas{\Delta x_1} \leq \beta_g \cdot u + O(u^2),\\
&= g(x+ \Delta x_1 + \Delta x_2), \qquad \meas{\Delta x_2} \leq \kappa_{g^{-1}} \cdot \meas{\Delta w} + O(\meas{\Delta w}^2).
\end{align*}
Since $\Delta x_1$ and $\Delta x_2$ are both perturbations of the same quantity $x$ there holds the triangle inequality for relative errors,
\[
\Delta x = \Delta x_1 + \Delta x_2,\qquad \meas{\Delta x} \leq \meas{\Delta x_1} + \meas{\Delta x_2} \leq (\beta_g + \kappa_{g^{-1}} \cdot \beta_h)\, u + O(u^2),
\]
and we get the assertion.
\end{proof}

Thus, we may complement the maxim from the beginning of this section by the following rule:
\medskip
\begin{quote}
\emph{If $g^{-1}$ is well-conditioned,
backward stable realizations of $g$ and $h$ induce a backward stable realization of $f$.}
\end{quote}
\medskip

Summarizing, the logical status of the proposed model is a follows. If the model predicts instability, we can expect
instability in reality---independently of how $g$ and $h$ are realized in practice.
Most probably, in examples that realize the worst case scenario of the condition number bound, there will be instability even if $g$
and $h$ were calculated exactly; a fact, which certainly shakes our faith in the algorithm. On the other hand, if the model predicts stability, the actual stability of
the algorithm depends on how $g$ and $h$ are  realized algorithmically. In the framework of backward stability,
stability of the realization of $g$ and $h$ implies stability of the resulting algorithm for $f$.


\paragraph{Example} Let us illustrate these points by reconsidering the example of {\S}\ref{sect2:examples}.A.
Here, we decompose $f(x) = \log^2(1+x)$, $x \approx 0$,  differently into
\[
f : x  \xmapsto{g} w = \log(1+x) \xmapsto{h} w^2.
\]
Now, the critical map $g^{-1} : w \longmapsto e^w-1$ has relative condition number $\kappa_{g^{-1}}\approx 1$ for $w \approx 0$.
The model alone would therefore predict numerical stability.
On the other hand, the full, recursive analysis has to take the actual algorithms for $g$ and $h$ into account. Step $h$, as a multiplication in IEEE arithmetic,  is certainly backward stable.
However, the status of $g$ is far less clear. If its realization is chosen to be based on the decomposition $g : x \mapsto z=1+x \mapsto \log(z)$,
then an analysis similar to {\S}\ref{sect2:examples}.A reveals instability. Otherwise, if $g$ is realized,
for instance, by using Kahan's stable algorithm as implemented in Matlab's {\tt log1p} command, the resulting algorithm for $f$ is stable, too.

Hence, the choice of the decomposition will critically determine the success or failure of the model. In general, making a
conclusive choice will depend on the user's experience or luck. However, we will show in the rest of the paper, that quite natural such decompositions
occur in the analysis of the stability of Gaussian elimination.

\section{The Backward Error of Linear Systems}\label{sect:oettli}

To prepare for a more {quantitative} analysis of algorithms for the solution of linear systems
of equations $Ax = b$ we recall the concept of the {backward error} of an output vector $\tilde{x} \in \R^m$.
{\em Normwise} analysis considers\footnote{Throughout the paper we deal with \emph{monotone} vector norms like the $1$-, $2$-, or $\infty$-norm,
and the induced matrix norms.}
\[
\eta = \min_{\,E \in \R^{m \times m}}\left\{\frac{\|E\|}{\|A\|} :  (A+E)\tilde{x} = b \,\right\},
\]
whereas {\em componentwise} analysis studies
\[
\omega = \min_{E \in \R^{m \times m}}\left\{\max_{ij}\frac{|E|_{ij}}{|A|_{ij}} : (A+E)\tilde{x} = b\,\right\}.
\]
The classic results \cite{Rigal} and \cite{Oettli} show that $\eta$ and $\omega$
can be calculated from the data of the linear system and the output vector $\tilde{x}$ by means
of the following simple formul{\ae}:
\begin{equation}\label{eq:wilk}
\eta = \frac{\|r\|}{\|A\|\cdot\|\tilde{x}\|},\qquad \omega = \max_{j=1:m}\frac{|r_j|}{(|A|\cdot |\tilde{x}|)_j}.
\end{equation}
Here, $r=b-A \tilde{x}$ denotes the \emph{residual} of $\tilde{x}$.
These formul{\ae}, which have very short and straightforward proofs \cite[pp.~120/122]{Higham}, are also valuable for the {\em a posteriori} assessment of computed solutions.
We will use them as a convenient point of departure for a quantitative analysis in the frame of our proposed model.

\section{Model Analysis of the Na\"{\i}ve Algorithm for Linear Systems}\label{sect:naiv}

As discussed in {\S}\ref{sect2:examples}, the na\"{\i}ve algorithm for the solution of a linear system
is given by the decomposition
\[
f:\;A  \xmapsto{\;\, g \;\,} B = A^{-1}  \xmapsto{\;\, h \;\,} x = B\cdot b.
\]
Our model analyzes how roundoff in $B$ affects the solution $x$ and its backward error:
\[
\tilde{f}:\; A  \xmapsto{\;\, g \;\,} B = A^{-1}  \xmapsto{\;\, fl \;\,} \tilde{B}=B+\Delta B \xmapsto{\;\, h \;\,} \tilde{x} = \tilde{B}\cdot b.
\]
The perturbation $|\Delta B| \leq u \cdot |B|$ induces,
by propagating
backwards through $g^{-1}$, an equivalent perturbation $\tilde{A} = A + \Delta A = g^{-1}(\tilde B)$ of the input matrix. By construction, we have $\tilde{A} \tilde{x} = b$,
\[
(A+\Delta A)(A^{-1} + \Delta B) = I,\qquad\text{i.e.,}\qquad \Delta A = -A \cdot \Delta B \cdot A - \Delta A\cdot \Delta B \cdot A,
\]
and therefore the componentwise estimate
\[
|\Delta A \cdot \tilde{x}| \leq |A|\cdot|A^{-1}|\cdot |A\tilde{x}|\, u + O(u^2).
\]
Since $r=b-A\tilde{x} = \Delta A \cdot \tilde{x}$ and $\tilde x = x + O(u)$, we get by (\ref{eq:wilk})
\begin{equation}\label{eq:bound1}
\eta = \frac{\|\Delta A \cdot \tilde{x}\|}{\|A\|\cdot\|\tilde{x}\|} \leq
\frac{\|\,|A|\cdot|A^{-1}|\cdot |Ax|\,\|}{\|A\|\cdot\|x\|}\,u + O(u^2) =:\gamma(A,x)\, u + O(u^2).
\end{equation}
To relate with better known quantities, we may further estimate
\[
\gamma(A,x) \;\leq  \;\|\,|A|\cdot|A^{-1}|\,\|=\cond(A^{-1}),
\]
in agreement with our qualitative analysis of {\S}\ref{sect2:examples}.B.
Thus, instability in the normwise concept appears to be only possible for badly conditioned matrices.

\subsection{Examples}\label{subsect.naiv}\hspace*{-0.2cm}\footnote{If not explicitly stated otherwise, all the examples in this paper use the norm $\|\cdot\|_\infty$.}
\paragraph{A} A notoriously badly conditioned matrix is the famous \emph{Hilbert matrix} $H_m$ for larger dimensions $m$.
In Matlab there is the command {\tt invhilb} that supplies $H_m^{-1}$ and allows to implement the na\"{\i}ve algorithm:\footnote{Here, and in the examples to follow,
we have cross-checked the \emph{actually calculated} backward errors with higher precision arithmetic. The first digits were always correct, so that the
conclusions we draw are not affected by roundoff errors in the computed residuals.}
\begin{verbatim}

>> m = 20; A = hilb(m); B = invhilb(m); b = ones(m,1); x = B*b;
>> eta = norm(b - A*x,inf)/norm(A,inf)/norm(x,inf)

eta = 1.2787e-005

\end{verbatim}
Thus, the na\"{\i}ve algorithm is unstable
as predicted by the a priori bound (\ref{eq:bound1}), which
turns out to be
\[
\eta = 1.27\cdots \times 10^{-5} \leq \gamma(A,x) \cdot u = 5.69\cdots \times 10^{-4};
\]
a fairly good prediction indeed. On the other hand, we have to be careful to base a prediction on coarser upper bounds
that were introduced for the ease of interpretation:
the condition number yields
\[
\eta \leq \cond(A^{-1}) \cdot u = 6.63\cdots \times 10^{11},
\]
which gives too pessimistic a picture of the actual backward error.

\paragraph{B} The following example \cite[p.~509]{Skeel1} shows that the na\"{\i}ve algorithm can be stable for \emph{some} badly conditioned matrices:
\[
A = \begin{pmatrix}
1 & 1 & -1 & -1 \\
1 &  0 &   0 &  1 \\
0  & 0  & \epsilon & 0\\
 0 & \epsilon & 0 & 0
\end{pmatrix},\qquad b = \begin{pmatrix}
0 \\
2 \\
1 \\
1\\
\end{pmatrix}, \qquad x = \begin{pmatrix}
1 \\
\epsilon^{-1} \\
\epsilon^{-1} \\
1\\
\end{pmatrix}.
\]
This matrix fulfills
\[
\cond(A) = 4,\qquad \cond(A^{-1}) = 2 + 4\epsilon^{-1}.
\]
However, numerical experiments with various small $0 < \epsilon \ll 1$ exhibit very small backward errors of about the size of the unit roundoff.
This is fully reflected by our model analysis, since
\[
\gamma(A,x) = 1 + \frac{\epsilon}{2} \approx 1.
\]

\section{Model Analysis of Gaussian Elimination: The Normwise Case}\label{sect:GEnorm}

As discussed in {\S}\ref{sect2:examples}.C the solution of a linear system $Ax=b$
by Gaussian Elimination corresponds to the decomposition
\[
f:\;A  \xmapsto{\;\, g \;\,} (L,U)  \xmapsto{\;\, h \;\,} x .
\]
In the model roundoff affects only the intermediate result, the $LU$-factorization, by
\[
\tilde f: \; A \xmapsto{g} (L,U) \xmapsto{fl} (\tilde L,\tilde U) = (L+\Delta L,U+\Delta U) \xmapsto{h} \tilde x.
\]
Here, the perturbations $|\Delta L | \leq u \cdot |L|$, $|\Delta U| \leq u \cdot |U|$
induce, by propagating through the inverse of $g$ (that is, matrix multiplication), an equivalent perturbation of the input matrix
\[
A + \Delta A =  \tilde{L} \tilde{U}  = (L+\Delta L)\cdot(U+\Delta U),\quad\text{i.e.,}\quad
\Delta A = \Delta L\cdot U + L\cdot \Delta U + \Delta L\cdot\Delta U.
\]
This way we obtain the componentwise estimate
\begin{equation}\label{eq:bound2}
|\Delta A| \leq   2 \,|L| \, |U|\cdot {u_*},\qquad {u_*} = u+u^2/2.
\end{equation}
Because of $r=b-A\tilde{x} = \Delta A \cdot \tilde{x}$ we get by (\ref{eq:wilk})
\begin{equation}\label{eq:bound3}
\eta  \leq \frac{\norm{\,|\Delta A|\cdot |\tilde{x}|\,}}{\norm{A}\cdot\norm{\tilde{x}}}
\leq 2\, \frac{\norm{\,|L|\cdot|U|\cdot|\tilde{x}|\,}}{\norm{A}\cdot\norm{\tilde{x}}} \,{u_*}
\leq 2\, \frac{\norm{\,|L|\cdot|U|\,}}{\|A\|} \,{u_*} =: 2 \,\gamma(L,U)\, {u_*},
\end{equation}
in agreement with our qualitative analysis of {\S}\ref{sect2:examples}.C.
If we restrict ourselves to monotone matrix norms, we can further estimate the \emph{growth factor} $\gamma(L,U)$ by using $U = L^{-1}\cdot A$
\[
\gamma(L,U) = \frac{\norm{\,|L|\cdot|U|\,}}{\|A\|} \leq \frac{\norm{\,|L|\cdot|L^{-1}|\,}\cdot\norm{A}}{\|A\|} = \cond(L^{-1}).
\]
Thus, an instability of Gaussian elimination in the normwise case requires a badly conditioned $L$-factor of the matrix $A$.
\subsection{Examples}
\paragraph{A} It is well known that Gaussian elimination without pivoting is bound to be \emph{unstable} for small pivot elements.
An example is given by
\[
A = \begin{pmatrix}
\epsilon & 1\\
1 & 1
\end{pmatrix},\qquad L =
\begin{pmatrix}
1 & 0\\
\epsilon^{-1} & 1
\end{pmatrix},\qquad
U  =
\begin{pmatrix}
\epsilon & 1\\
0 & 1 - \epsilon^{-1}
\end{pmatrix}.
\]
For $\epsilon = u$, $b= (1,0)^T$ a numerical experiment yields\footnote{We write $a\doteq b$
if $a-b\approx u$.} $\tilde{x} \doteq (-2,1)$;  the exact solution, however, would be $x  \doteq (-1,1)^T$.
The backward error turns out to be $\eta \doteq \frac{1}{4}$. On the other hand we have
\[
\gamma(L,U) = \epsilon^{-1},\qquad \cond(L^{-1}) = 1+2\epsilon^{-1},
\]
which,  by (\ref{eq:bound3}), gives the fairly good prediction $\eta \leq \epsilon^{-1}\cdot u = 1$.

\paragraph{B} Gaussian elimination with partial pivoting yields an $L$-factor that satisfies $|L|\leq 1$ componentwise.
This can be used \cite[p.~143]{Higham} to show that
\[
\gamma(L,U) \leq \cond(L^{-1}) \leq 2^{m}-1,
\]
which proves that the growth factor remains bounded for \emph{fixed} dimension $m$. However, the upper bound on $\cond(L^{-1})$ is attained for
Wilkinson's famous matrix
\[
A = \begin{pmatrix}
1 & & &1\\
-1 & \ddots & &\vdots\\
\vdots & \ddots & \ddots & \vdots\\
-1 & \cdots & -1 & 1
\end{pmatrix},\qquad L = \begin{pmatrix}
1 & & \\
-1 & \ddots & \\
\vdots & \ddots & \ddots \\
-1 & \cdots & -1 & 1
\end{pmatrix}.
\]
Numerical experiments quickly exhibit very large backward errors:
\begin{verbatim}

>> m = 53; A = eye(m)-tril(ones(m),-1); A(:,m) = 1;
>> rand('seed',42); b = rand(m,1); x = A\b;
>> eta = norm(b-A*x,inf)/norm(A,inf)/norm(x,inf)

eta = 3.2342e-003

\end{verbatim}
Our analysis yields a fairly good prediction,
\[
 \eta = 3.23\cdots\times 10^{-3} \leq 2 \,\gamma(L,U)\cdot {u_*} = 3.77\cdots \times 10^{-2}.
\]

\paragraph{C} For symmetric positive definite matrices, the solution of the linear system $Ax=b$
by \emph{Cholesky factorization} corresponds to the decomposition
\[
f:\;A  \xmapsto{\;\, g \;\,} L  \xmapsto{\;\, h \;\,} x
\]
with $A = L\cdot L^T$.
A perturbation $\tilde{L} = L + \Delta L$ of the intermediate result by roundoff,
\[
|\Delta L | \leq u \cdot |L|,
\]
induces, as for (\ref{eq:bound3}), the backward error (with resprect to the norm $\norm{\cdot}_2$)
\[
\eta \leq 2\,\frac{\norm{\,|L|\cdot|L^T|\,}_2}{\norm{A}_2}\, {u_*} = 2\, \gamma(L,L^T)\,{u_*}.
\]
Since $\norm{\,|L|\,}_2 \leq \sqrt{m}\|L\|_2$ for any $m\times m$ matrix, we infer
 \cite[p.~198]{Higham}
\[
\gamma(L,L^T) \leq \frac{\norm{\,|L|\,}_2\norm{\,|L^T|\,}_2}{\norm{A}_2} \leq m \,\frac{\norm{L}_2\norm{L^T}_2}{\norm{LL^T}_2} = m.
\]
Hence, we have
\[
\eta \leq 2\,m\,{u_*},
\]
which hints to the perfect normwise backward stability of the Cholesky method.

\section{Model Analysis of Gaussian Elimination: The Componentwise Case}\label{sect:GEcomp}

The matrix estimate (\ref{eq:bound2}) immediately yields an estimate of the \emph{componentwise} backward error,
\begin{multline}\label{eq:chain}
\omega \;=\; \max_j \frac{|\Delta A\cdot \tilde{x}|_j}{(|A|\cdot|\tilde{x}|)_j} \;\leq\; \max_j \frac{(|\Delta A|\cdot |\tilde{x}|)_j}{(|A|\cdot|\tilde{x}|)_j}
\;\leq\; 2\, \max_j \frac{(|L|\cdot|U|\cdot|\tilde{x}|)_j}{(|A|\cdot|\tilde{x}|)_j}\, {u_*}\\*[2mm]
\;\leq \;2\, \frac{\max_j (|L|\cdot|U|\cdot|\tilde{x}|)_j}{\min_j (|A|\cdot|\tilde{x}|)_j)} \;{u_*} \;=\;
2\, \frac{\norm{\,|L|\cdot|U|\cdot|\tilde{x}|\,}_\infty}{\norm{\,|A|\cdot|\tilde{x}|\,}_\infty} \underbrace{\frac{\max_j (|A|\cdot|\tilde{x}|)_j}{\min_j (|A|\cdot|\tilde{x}|)_j}}_{=\sigma(A,\tilde{x})} \;{u_*},
\end{multline}
which by $U = L^{-1}A$, that is $|U| \leq |L^{-1}|\cdot |A|$, induces \cite[Thm.~4.4]{Skeel1}
\begin{equation}\label{eq:bound4}
\omega \leq 2\,\cond(L^{-1}) \, \sigma(A,\tilde{x})\, {u_*}.
\end{equation}
As our derivation shows, this is not necessarily the best possible concise bound, but it allows for the easy comparison with the normwise bound (with respect to $\norm{\cdot}_\infty$)
\[
\eta \leq 2 \, \cond(L^{-1}) \,  {u_*}.
\]
We see that the componentwise bound just differs by the additional factor $\sigma(A,\tilde{x})\geq 1$. This factor measures the quality of
the \emph{scaling} of the linear system with respect to $\tilde{x}$ and predicts an instability for badly scaled systems.

\subsection{Examples}
\paragraph{A} We return to the example of {\S}\ref{subsect.naiv}.B. The growth factor and the scaling are given
by
\[
\cond(L^{-1}) = 3+4\epsilon,\qquad \sigma(A,x) = 2 + 2 \epsilon^{-1}.
\]
Experimentally, for $\epsilon = 10^{-16}$, Gaussian elimination yields (partial pivoting is not used here because of $|L|\leq 1$)
\[
\eta = 2.84\cdots \times 10^{-17} \leq 2\,\cond(L^{-1}) {u_*} = 6.66\cdots \times 10^{-16}.
\]
On the other hand, the componentwise backward error satisfies
\[
\omega = 0.499\cdots \leq 2\,\cond(L^{-1}) \sigma(A,x){u_*} = 13.3\cdots.
\]
Thus, the model analysis helps to understand the actual behavior of the two error concepts. In particular, we see that scaling
can be an issue for Gaussian elimination with partial pivoting if analyzed componentwise.

\paragraph{B} There are matrices, for which the upper bound (\ref{eq:bound4}) turns out to be too coarse. As an example, we consider
totally positive matrices $A$ such as the Hilbert matrix of {\S}\ref{subsect.naiv}.A or matrices that appear in spline interpolation. These matrices factor with $L \geq 0$ and $U\geq 0$. Thus, we best stay with the
following intermediate step in the chain of estimates (\ref{eq:chain}):
\[
\omega \;\leq\; 2 \max_j \frac{(|L|\cdot|U|\cdot|\tilde{x}|)_j}{(|A|\cdot|\tilde{x}|)_j}\, {u_*}.
\]
Here, we obviously have  $|L|\cdot |U| = |A|$ and we can therefore directly infer the perfect stability estimate \cite{deBoor}
\[
\omega \;\leq\; 2\,{u_*}.
\]

\section{Model Analysis of a Single Iterative Refinement Step}\label{sect:GESkeel}

In this final section we will apply the model analysis to the understanding of the
results \cite{Skeel2} on iterative refinement of Gaussian elimination. We recall that the iterative refinement
of a calculated solution $\tilde x$ of a linear system $Ax=b$ consists of three steps: compute the residual
$r_0 = b -A\tilde x$, solve $Aw=r_0$ for a calculated correction $\tilde w$ (reusing the $LU$-decomposition of $A$), update $\tilde y=\tilde x + \tilde w$. If there were
no roundoff errors in the refinement steps (that is, $\tilde w = w$), we would obtain $\tilde y=x$, the exact solution.

In the previous two sections, the model analysis of Gaussian elimination allowed for roundoff errors just in the $L$- and $U$-factors of $A$
yielding some equivalent perturbation of that matrix.
Because of the reuse of these factors in the iterative refinement step, we reasonably assume that both Gaussian elimination steps, that is,
those leading to
 $\tilde x$ and  $\tilde w$, are affected by roundoff through
a \emph{single} perturbation  $\tilde{A} = A+ \Delta A$ satisfying the estimate (\ref{eq:bound2}). This way, the result $\tilde y$ of the
 iterative refinement  is given by
\[
\tilde{y} = \tilde{x} + \tilde w,\qquad
(A+\Delta A)\tilde w = r_0 = b- A \tilde{x} = \Delta A \cdot \tilde{x}.
\]
The residual after this step is $r_1 = b - A\tilde{y} = r_0- A \tilde w = \Delta A \tilde w$,
and therefore
\[
A \tilde w = r_0- \Delta A  \tilde w = \Delta A(\tilde{x}-\tilde w) = \Delta A (\tilde{y}-2\tilde w),\qquad \tilde w = A^{-1} \Delta A \tilde{y} -2A^{-1}
\Delta A \tilde w.
\]
Hence we have
\[
|\Delta A \tilde w| \;\leq\; |\Delta A|\, |A^{-1}|\,|\Delta A| \,|\tilde{y}|\; + \;2\, |\Delta A|\,|A^{-1}|\,|\Delta A \tilde w|,
\]
which by (\ref{eq:bound2}), that is $|\Delta A| \leq 2\,|L|\,|U|\,{u_*}
\leq 2\,|L|\,|L^{-1}|\,|A|\,{u_*}$,
implies
\begin{multline*}
\|\Delta A\,\tilde w\|_\infty \;\leq\; 4\,\cond^2(L^{-1})\,\cond(A^{-1})\cdot
\norm{\,|A|\,|\tilde{y}|\,}_\infty\,u_*^2\\
+\; 4\,\cond(L^{-1})\,\cond(A^{-1})  \,{u_*}\, \norm{\Delta A\,\tilde w}_\infty.
\end{multline*}
If $4\,\cond(L^{-1})\,\cond(A^{-1})  \,{u_*} < 1$ we can solve for $\|\Delta A\,w\|_\infty$ and get---as in the derivation
of (\ref{eq:chain})---the following upper bound of the backward error of $\tilde{y}$:
\begin{multline}\label{eq:boundSkeel}
\omega_1 \;=\; \max_j \frac{|r_1|_j}{(|A|\cdot|\tilde{y}|)_j} \;=\; \max_j \frac{|\Delta A\cdot \tilde w|_j}{(|A|\cdot|\tilde{y}|)_j} \;\leq\; \frac{\max_j |\Delta A\cdot \tilde w|_j}{\min_j (|A|\cdot|\tilde{y}|)_j}\\*[2mm]
=\;\frac{\norm{\Delta A \,\tilde w}_\infty}{\norm{\,|A|\,|\tilde{y}|\,}_\infty} \sigma(A,\tilde{y})
\leq \;\frac{4\,\cond^2(L^{-1})\,\cond(A^{-1})\,\sigma(A,\tilde{y})\, {u_*}}{1-4\,\cond(L^{-1})\,\cond(A^{-1})  \,{u_*}}\; {u_*}.
\end{multline}
Because of $\cond(L^{-1}) \geq 1$, $\sigma(A,\tilde{y})\geq 1$, the premise
is in particular satisfied if
\begin{equation}\label{eq:condition}
8 \,\cond^2(L^{-1})\,\cond(A^{-1})\,\sigma(A,\tilde{y}) \,{u_*} \leq 1,
\end{equation}
for which we obtain from (\ref{eq:boundSkeel}) the simple perfect bound $\omega_1 \leq {u_*}$.
Except for a constant depending on the dimension $m$ this
is \emph{exactly} the result \cite[p.~239]{Higham} of an elaborate analysis that takes all the details of
roundoff error rigorously into account.\footnote{The catch, of course, is that without doing the full analysis we would not
know if we had really determined the full bound. However, the point of this paper is a better understanding of the underlying mathematical
structure. If, by neglecting many details, we come to predict the same bounds with much less effort we seem to have put the focus
on the right spot.}

Summarizing our analysis predicts: As long as the linear system is not too badly conditioned ($\cond(A^{-1})$ is not too large)
and not too badly scaled ($\sigma(A,\tilde{y})$ is not too large), and
Gaussian elimination is not too unstable ($\cond(L^{-1})$ is not too large), \emph{one step of iterative refinement implies componentwise
backward stability}.

\begin{figure}[tbp]
\begin{center}
\includegraphics[scale=0.385,angle=270]{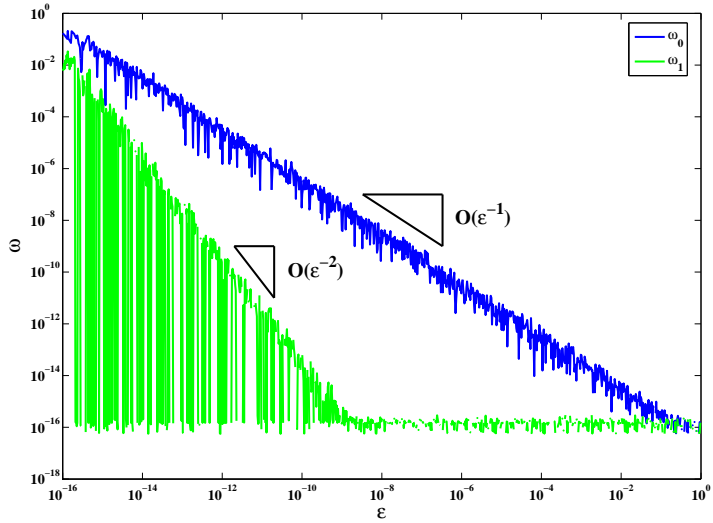}
\caption{Backward errors $\omega_0$ and $\omega_1$ vs. $\epsilon$}\label{fig.1}
\end{center}
\end{figure}

\subsection{Example}
We consider the example \cite[p.~500]{Skeel1}
\[
A = \begin{pmatrix}
3 & 2 & 1\\
2 & 2\epsilon & 2\epsilon \\
1 & 2\epsilon & -\epsilon
\end{pmatrix},
\qquad
b = \begin{pmatrix}
3 + 3\epsilon\\
6\epsilon\\
2\epsilon
\end{pmatrix},
\qquad
x = \begin{pmatrix}
\epsilon\\
1\\
1
\end{pmatrix},
\]
of  a  well conditioned (for this particular right hand side $b$), but badly scaled linear system. Because of
\[
\cond(A^{-1}) = \frac{6}{5}\epsilon^{-1} + O(1),\quad \sigma(A,x) = \frac{3}{4}\epsilon^{-1}+O(1),\quad \cond(L^{-1}) = \frac{8}{3} + O(\epsilon),
\]
 condition (\ref{eq:condition}) reads as
\[
1\geq 8 \,\cond^2(L^{-1})\,\cond(A^{-1})\,\sigma(A,\tilde{y}) \,{u_*}
= \frac{256}{5}\epsilon^{-2} +O(\epsilon^{-1}),
\]
that is, one step of iterative refinement is predicted to imply stability as long as $\epsilon$
remains  larger than about the square root of the unit roundoff,
\[
\epsilon \geq \frac{16}{5} \sqrt{5u} + O(u^2) \approx 7.5\times 10^{-8}.
\]
In fact, the upper bound (\ref{eq:bound4}) predicts that the componentwise backward error $\omega_0$
of $\tilde{x}$ behaves like $\omega_0 = O(\epsilon^{-1}u)$; whereas the upper bound (\ref{eq:boundSkeel})
predicts $\omega_1 = O(\epsilon^{-2}u^2)$ for the first refinement
step $\tilde{y}$. All this can perfectly be observed in an actual numerical experiment, see Figure~\ref{fig.1}.

\medskip

\paragraph{Acknowledgements}
We are grateful to Nick Higham for commenting on a draft of this manuscript.
\bibliographystyle{kluwer}
\bibliography{stability}

\end{document}